\documentclass[11pt]{article}
\usepackage{amsmath}
\usepackage{amssymb}
\usepackage{url}
\usepackage[pstarrows]{pict2e}
\newcommand\dstyle\displaystyle
\newcommand\nonu{\nonumber}
\newcommand\sa{\smallskipamount}
\newcommand\ma{\medskipamount}
\newcommand\ba{\bigskipamount}
\newcommand\sLP{\\[\sa]}
\newcommand\mLP{\\[\ma]}
\newcommand\bLP{\\[\ba]}
\newcommand\ga\gamma
\newcommand\la\lambda

\newcommand\bPP{\\[\ba]\indent}

\newcommand\CC{\mathbb{C}}
\newcommand\RR{\mathbb{R}}
\newcommand\ZZ{\mathbb{Z}}
\newcommand\al\alpha
\newcommand\be\beta
\newcommand\de\delta
\newcommand\si\sigma

\newcommand\wt{\widetilde}
\newcommand\half{\frac12}
\newcommand\thalf{\tfrac12}
\newcommand{\hyp}[5]{\,\mbox{}_{#1}F_{#2}\!\left(
  \genfrac{}{}{0pt}{}{#3}{#4};#5\right)}
\newcommand{\qhyp}[5]{\,\mbox{}_{#1}\phi_{#2}\!\left(
  \genfrac{}{}{0pt}{}{#3}{#4};#5\right)}
\newcommand{\eup}{\textup{e}}
\newcommand{\iup}{\mkern1mu\textup{i}\mkern1mu}
\newcommand\RHS{right-hand side}

\numberwithin{equation}{section}
\newtheorem{theorem}{Theorem}[section]
\newtheorem{Remark}[theorem]{Remark}
\newenvironment{remark}{\begin{Remark}\rm}{\end{Remark}}
\newtheorem{Example}[theorem]{Example}

\begin{document}
\title{Charting the $q$-Askey scheme. III. Verde-Star scheme
for $q=1$}
\author{Tom H. Koornwinder}
\date{Dedicated to the memory of Jos\'e Carlos Petronilho,\\
who loved classical orthogonal polynomials}
\maketitle

\begin{abstract}
Following Verde-Star, Linear Algebra Appl. 627 (2021), we label families
of orthogonal polynomials in the $q=1$ Askey scheme together with their
hypergeometric representations by three sequences $x_k, h_k, g_k$
of  polynomials in $k$,
two of degree 2 and one of degree 4, satisfying certain constraints.
Except for the Hermite polynomials, this gives rise to a precise classification and
a very simple uniform parametrization of these
families together with their limit transitions. This is displayed in a graphical
scheme.
We also discuss limits from the $q$-case to the case
$q=1$, although this cannot be done in a uniform way.
\end{abstract}

\section{Introduction}
This is the third one in a series of papers which aim to reparametrize the
families of orthogonal polynomials in the
$q$-Askey scheme \cite[Chapter 14]{2} and its limit cases
for $q\to1$ \cite[Appendix]{6}, \cite[Chapter 9]{2} and for $q\to-1$ \cite{16},
\cite{27} in a uniform way
following Verde-Star's set-up \cite{1} (see also Vinet \& Zhedanov \cite{9}),
and to redraw the ($q$-)Askey scheme
in accordance with these uniform parameters. Beside parameters from \cite{1}
we also take parameters from the structure constants in the Zhedanov algebras
\cite{4}
corresponding to the families in the ($q$-)Askey scheme.
Our first paper \cite{25} applied \cite{1} to the general $q$-case, yielding the
so-called $q$-Verde-Star scheme.
Our second paper \cite{26}, also for general $q$,
related the Verde-Star parameters to the
structure constants of the Zhedanov algebra, thus yielding the so-called
$q$-Zhedanov scheme. The present paper considers the Verde-Star parameters in
the $q=1$ case, thus yielding the so-called Verde-Star scheme.

Verde-Star's approach uses data $h_k$, $x_k$ and $g_k$ which are, for $q=1$,
polynomials in $k$ of maximal degrees 2, 2 and 4, respectively. The $x_k$
are nodes of monic Newton polynomials $v_n(x)$. A linear operator $L$ acting
on the space of polynomials is defined by requiring that
$Lv_n=h_n v_n+g_n v_{n-1}$. Then define monic polynomial $u_n(x)$ of degree $n$
such that $Lu_n=h_n u_n$. Polynomials $u_n$ are next classified such that they
satisfy a three-term recurrence relation. Families are characterized by
the degrees in $k$ of the corresponding data $h_k,x_k,g_k$. This approach
disregards orthogonality properties, so it is different from the usual
approach in Bochner type classifications.
Also observe that one family obtained in the present approach can stand
for several families of orthogonal polynomials,
for which sets on which they are orthogonal differ.
But very prominent in this approach is a hypergeometric expression written in
a conceptual way. Indeed, up to a constant factor $u_n(x)$ equals
\[
\sum_{k=0}^n\,\prod_{j=0}^{k-1}\frac{(h_n-h_j)(x-x_j)}{g_{j+1}}\,.
\]

As already observed in \cite{1}, the Hermite polynomials do not fit into this set-up.
As could be expected, the Verde-Star scheme is much less complicated than its
$q$-analogue \cite{25}.
But it is very satisfactory that, just as in the $q$-case and apart from the case of
the Hermite polynomials, the
Verde-Star scheme confirms in a conceptual way that the Askey scheme
\cite[Chapter 9]{2} precisely gives all families and all limit transitions
which lower the number of parameters by one. Nicer also than in the $q$-case,
we can describe all families in the Verde-Star scheme in a uniform way with
four parameters.
We did not succeed to give a uniform parametrization for the general $q$-case
such that also limits for $q\to1$ can be taken. Such limits have to be made
in an ad hoc way.

The contents of this paper are as follows.
In Section 2 we reacapitulate Verde-Star's approach, both in the general $q$-case
and for $q=1$, and we illustrate it with the cases of Askey--Wilson polynomials
and Wilson polynomials. In Section 3, which is the heart of this paper, we
present the Verde-Star scheme and we give the uniform parametrization.
In Section 4 we list limits from families in the $q$-Verde-Star scheme
to corresponding families being on the same row in the Verde-Star scheme.
Finally, in Appendix A we give explicit data for all families in the Verde-Star
scheme and in Appendix B we give some worked-out examples of limits for $q\to1$.

\paragraph{Note}
For definition and notation of ($q$-)shifted factorials and
($q$-)hypergeometric series we follow \cite[\S1.2]{21}.
We will only need terminating series.
For formulas on orthogonal polynomials in the ($q$-)Askey scheme we
refer to \cite[Chapters 9 and 14]{2}.

\section{Verde-Star's theorem}
In this section we summarize Verde-Star's results \cite{1}, and we will make them
concrete in the case of Askey--Wilson polynomials and of Wilson polynomials.
The case $q=-1$ will not be considered here.
\subsection{The general set-up and two examples}
Let $h_k$, $x_k$, $g_k$ ($k=0,1,2,\ldots$) be given sequences such that all $h_k$
are distinct. Define monic Newton type polynomials
\begin{equation}
v_k(x):=(x-x_0)(x-x_1)\ldots(x-x_{k-1})\;\;(k\ge1),\quad v_0(x)=1,
\label{62}
\end{equation}
(here it is allowed that $x_k=c$ for all $k$, by which $v_k(x)=(x-c)^k$),
and define monic polynomials
\begin{equation}
u_n(x):=\sum_{k=0}^n c_{n,k}\,v_k(x),
\label{63}
\end{equation}
where
\begin{equation}
c_{n,k}=\prod_{j=k}^{n-1}\frac{g_{j+1}}{h_n-h_j}\;(0\le k<n),\quad
c_{n,n}=1.
\label{64}
\end{equation}
Define a linear operator $L$ on the space of polynomials by
\begin{equation}
Lu_n=h_nu_n\quad(n\ge0).
\label{65}
\end{equation}
Then
\begin{equation}
Lv_n=h_nv_n+g_nv_{n-1}\;(n>0),\quad Lv_0=h_0v_0.
\label{66}
\end{equation}
Note that any two of equalities \eqref{64}, \eqref{65}, \eqref{66} implies the
third equality.

\begin{Example}[Askey--Wilson polynomials]\rm
Let $q\ne0$, $1\notin q^\ZZ$, $a\ne0$, $1\notin abcdq^{\ZZ_{\ge0}}$. Let
\begin{align}
u_n(x)&=
\frac{(ab,ac,ad;q)_n}{a^n\,(q^{n-1}abcd;q)_n}\,
\qhyp43{q^{-n},q^{n-1}abcd,az,az^{-1}}{ab,ac,ad}{q,q}\quad(x=z+z^{-1})
\label{67}\\*
&=\frac1{(q^{n-1}abcd;q)_n}\,p_n\big(\thalf x;a,b,c,d\,|\,q\big),
\label{68}
\end{align}
where $p_n\big(x;a,b,c,d\,|\,q\big)$ is the Askey--Wilson polynomial in
usual notation and standardization \cite[(14.1.1)]{2}.
Then \eqref{67} can be rewritten as \eqref{63} with $v_k(x)$ and $c_{n,k}$
given by \eqref{62} and \eqref{64} and $h_k,x_k,g_k$ given by
\begin{equation}
\begin{split}
&h_k=q^{-k}(1-q^k)(1-abcdq^{k-1}),\quad
x_k=aq^k+a^{-1}q^{-k},\\
&g_k=a^{-1}q^{-2k+1}(1-abq^{k-1})(1-acq^{k-1})(1-adq^{k-1})(1-q^k).
\end{split}
\label{69}
\end{equation}
Then
\begin{equation}
v_k(x)=(-1)^k a^{-k} q^{-\half k(k-1)} (az,az^{-1};q)_k.
\label{73}
\end{equation}
With $h_k$ as in \eqref{69}, we can match \eqref{65} with
the eigenvalue equation \cite[(14.1.7)]{2}. So \eqref{66}
with \eqref{69} and \eqref{73} gives an alternative way to define
the second order $q$-difference operator in that eigenvalue equation.
\end{Example}

\begin{Example}[Wilson polynomials]\rm
Let $a+b+c+d\notin\ZZ_{\le0}$. For $x=y^2$ let
\begin{align}
u_n(x)
&=\frac{(-1)^n(a+b,a+c,a+d)_n}{(n+a+b+c+d-1)_n}\,
\hyp43{-n,n+a+b+c+d-1,a+\iup y,a-\iup y}{a+b,a+c,a+d}1
\label{70}\\
&=\frac{(-1)^n}{(n+a+b+c+d-1)_n}\,W_n(x;a,b,c,d),
\label{71}
\end{align}
where $W_n(x;a,b,c,d)$ is the Wilson polynomial in usual
notation and standardization \cite[(9.1.1)]{2}.
Then \eqref{70} can be rewritten as \eqref{63} with $v_k(x)$ and $c_{n,k}$
given by \eqref{62} and \eqref{64} and $h_k,x_k,g_k$ given by
\begin{equation}
\begin{split}
&h_k=-k(k+a+b+c+d-1),\quad x_k=-(a+k)^2,\\
&g_k=k(k+a+b-1)(k+a+c-1)(k+a+d-1).
\end{split}
\label{72}
\end{equation}
Then
\begin{equation}
v_k(x)=(a+\iup y,a-\iup y)_k.
\label{74}
\end{equation}
With $h_k$ as in \eqref{72}, we can match \eqref{65} with
the eigenvalue equation \cite[(9.1.6)]{2}. So also here \eqref{66}
with \eqref{72} and \eqref{74} gives an alternative way to define
the second order difference operator in that eigenvalue equation.
\end{Example}
\subsection{Verde-Star's theorem and some symmetries}
\label{91}
Following \cite{1} we now take some $\xi\ne0$
and assume that $g_k,x_k,h_k$ are solutions of difference equation of order
three, three and five, respectively, of the form
\begin{align}
h_{k+3}-h_k&=\xi(h_{k+2}-h_{k+1}),\label{75}\\
x_{k+3}-x_k&=\xi(x_{k+2}-x_{k+1}),\label{76}\\
g_{k+5}-g_k&=(\xi^2-\xi-1)\big(g_{k+4}-g_{k+1}-(\xi-1)(g_{k+3}-g_{k+2})\big),\quad
g_0=0.
\label{77}
\end{align}
Solutions must hold for $k\in\ZZ$, but they will be used only for $k\ge0$.
Solutions will be completely determined by the 10 values
$h_0,h_1,h_2$, $x_0,x_1,x_2$, $g_1,g_2,g_3,g_4$.

\begin{theorem}[Verde-Star {\cite[Theorem 6.1]{1}}]\label{79}
Let the sequences $h_k,x_k,g_k$ satisfy \eqref{75}, \eqref{76}, \eqref{77}. Then
$u_n(x)$ satisfies a three-term recurrence relation
\begin{equation}
xu_n(x)=u_{n+1}(x)+A_n u_n(x)+B_n u_{n-1}(x)\quad(n\ge1)
\label{80}
\end{equation}
if and only if
\begin{align}
g_3&=\xi\big((x_2-x_0)(h_2-h_0)-(x_1-x_0)(h_2-h_0)-(x_2-x_0)(h_1-h_0)+g_2-g_1\big),
\label{81}\\
g_4&=\xi(\xi-1)\Big(
(\xi+1)\big((x_1-x_0)(h_1-h_0)-(x_1-x_0)(h_2-h_0)-(x_2-x_0)(h_1-h_0)\big)
\nonu\\
&\qquad\qquad\qquad\qquad\qquad\qquad\qquad\qquad\qquad
+\xi(x_2-x_0)(h_2-h_0)+g_2\Big)+(1-\xi^2)g_1.\label{82}
\end{align}
\end{theorem}

\begin{remark}
Put $g_k=x_{k-1}(h_k-h_0)+d_k$. Suppose that \eqref{75} and \eqref{76} hold. Then
\eqref{77} holds if and only if $d_k$ satisfies
$d_{k+3}-d_k=\xi(d_{k+2}-d_{k+1})$ and $d_0=0$ (see \cite[\S10, first step]{1}).
\end{remark}

Because of the constraints \eqref{80} and \eqref{81},
solutions of \eqref{75}, \eqref{76}, \eqref{77}
will be completely determined by the 8 values
$h_0,h_1,h_2$, $x_0,x_1,x_2$, $g_1,g_2$.\\
Note that the right-hand sides of \eqref{81} and \eqref{82} are:
\begin{itemize}
\item
homogeneous of degree 1 in the $h_i$ and $g_i$;
\item
invariant under translations of all $h_i$ with the same amount;
\item
homogeneous of degree 1 in the $x_i$ and $g_i$;
\item
invariant under translations of all $x_i$ with the same amount;
\item
invariant under exchange of $x_i$ and $h_i$.
\end{itemize}
Corresponding to the first four above items we observe dilation and translation
symmetries in \eqref{62}--\eqref{66}:
\begin{enumerate}
\item
$h_k\to\rho h_k$, $g_k\to\rho g_k$ ($\rho\ne0$), which implies $L\to\rho L$.
\item
$h_k\to h_k+\si$, which implies $L\to L+\si$.
\item
$x_k\to\rho x_k$, $g_k\to\rho g_k$ ($\rho\ne0$), which implies\\
$v_k(x)\to \rho^k v_k(\rho^{-1} x)$,
$u_n(x)\to \rho^n u_n(\rho^{-1} x)$,
$L\to\de_\rho\circ L\circ\de_\rho^{-1}$, where
$(\de_\rho p)(x):=p(\rho^{-1}x)$.
\item
$x_k\to x_k+\si$, which implies $v_k(x)\to v_k(x-\si)$, $u_n\to u_n(x-\si)$,
$L\to\tau_\si\circ L\circ\tau_\si^{-1}$, where $(\tau_\si p)(x)=p(x-\si)$.
\end{enumerate}
We do not consider the application of these symmetries as essential changes.
So only four of the eight values $h_0,h_1,h_2$, $x_0,x_1,x_2$, $g_1,g_2$
are essential.

\begin{remark}
By \cite[(5.5), (5.6)]{1} the coefficients $A_n$ and $B_n$ in
the three-term recurrence relation \eqref{80}
are given by
\begin{align}
A_n&=x_n+\frac{g_{n+1}}{h_n-h_{n+1}}-\frac{g_n}{h_{n-1}-h_n}\,,
\label{14}\\
B_n&=\frac{g_n}{h_{n-1}-h_n}\left(\frac{g_{n-1}}{h_{n-2}-h_n}
-\frac{g_n}{h_{n-1}-h_n}+\frac{g_{n+1}}{h_{n-1}-h_{n+1}}+x_n-x_{n-1}\right).
\label{15}
\end{align}
For $n=0$ \eqref{80} and \eqref{14} degenerate to
\begin{equation*}
u_1(x)=x-A_0,\qquad A_0=x_0-\frac{g_1}{h_1-h_0}\,.
\end{equation*}
\end{remark}

\begin{remark}
In Theorem \ref{79} it is allowed that $B_n=0$ for all $n$.
This degenerate case
will certainly happen if $g_n=0$ for all~$n$.
Then $A_n=x_n$ and $u_n=v_n$, clearly
not belonging to a family of orthogonal polynomials. From now on we will exclude
this case.

It is also possible that the $B_n$ are zero because the second factor on
the \RHS\ of \eqref{15} is zero. We will allow this case.
\end{remark}

\begin{remark}
We may have that $g_n$ vanishes only for some values of $n$.
Let then $n=N+1$ be the lowest value of $n$ for which $g_n=0$. Then
$c_{n,k}=0$ if $N<k<n$. If we only consider $u_n$ for $n\le N$ we obtain
one of the finite systems of orthogonal polynomials in the ($q$-)Askey scheme.
\end{remark}

\subsection{$x\leftrightarrow h$ duality and ($q$-)hypergeometric form}
We already observed that the right-hand sides of \eqref{81} and \eqref{82} are
invariant under exchange of $x_i$ and $h_i$.
If all $x_i$ are distinct then the exchange of $x_i$ and $h_i$
relates $u_n$ given by \eqref{63} to its \emph{dual} $\wt u_n$
given by
\begin{equation}
\wt u_n(x)=\sum_{k=0}^n \wt c_{n,k}\,\wt v_k(x),\quad
\wt v_k(x)=\prod_{j=0}^{k-1}(x-h_j),\quad
\wt c_{n,k}=\prod_{j=k}^{n-1}\frac{g_{j+1}}{x_n-x_j}\,.
\label{83}
\end{equation}
If we put
\begin{align}
U_n(x)&=u_n(x)\,\prod_{j=0}^{n-1} \frac{h_n-h_j}{g_{j+1}}=
\sum_{k=0}^n\,\prod_{j=0}^{k-1}\frac{(h_n-h_j)(x-x_j)}{g_{j+1}}\,,\label{84}\\
\wt U_m(x)&=\wt u_m(x)\,\prod_{j=0}^{m-1} \frac{x_m-x_j}{g_{j+1}}=
\sum_{k=0}^m\,\prod_{j=0}^{k-1}\frac{(x_m-xj)(x-h_j)}{g_{j+1}}\label{85}
\end{align}
then
\begin{equation}
U_n(x_m)=\wt U_m(h_n).
\label{86}
\end{equation}

In general, while $u_n(x)$ is the monic polynomial, $U_n(x)$ is the
($q$-)hypergeometric form of that polynomial. For instance, if $h_n,x_n,g_n$ are
the data \eqref{69} for the Askey--Wilson polynomial \eqref{67}, \eqref{68}
then \eqref{84} becomes
\begin{align}\label{92}
U_n(x)=U_n(x;a,b,c,d;q)&=
\qhyp43{q^{-n},q^{n-1}abcd,az,az^{-1}}{ab,ac,ad}{q,q}\quad(x=z+z^{-1})\\*
&=\frac{a^n}{(ab,ac,ad;q)_n}\,p_n\big(\thalf x;a,b,c,d\,|\,q\big).\nonu
\end{align}
Also, if $h_n,x_n,g_n$ are
the data \eqref{72} for the Wilson polynomial \eqref{70}, \eqref{71}
then \eqref{84} becomes
\begin{align}\label{93}
U_n(x)=U_n(x;a,b,c,d)&=
\hyp43{-n,n+a+b+c+d-1,a+\iup y,a-\iup y}{a+b,a+c,a+d}1\quad(x=y^2)\\*
&=\frac1{(a+b,a+c,a+d)_n}\,W_n(x;a,b,c,d).\nonu
\end{align}

\subsection{The general $q$-case}
Put $\xi=1+q+q^{-1}$ with $q\ne0$ and assume $1\notin q^\ZZ$.
For this case see \cite[\S7]{1}.
The solutions of \eqref{75}, \eqref{76}, \eqref{77} take the form
\begin{equation}
\begin{split}
&h_k=a_{-1} q^{-k}+a_0+a_1 q^k,\qquad
x_k=b_{-1} q^{-k}+b_0+b_1 q^k,\\
&g_k=d_{-2} q^{-2k}+d_{-1} q^{-k}+d_0+d_1 q^k+d_2 q^{2k},\\
&a_{-1}\notin a_1 q^{\ZZ_{>0}},\quad
\sum_{i=-2}^2 d_i=0,\mbox{ $d_i\ne0$ for some $i$}.
\end{split}
\label{87}
\end{equation}
Furthermore, the constraints \eqref{81}, \eqref{82} now simplify to
\begin{equation}
d_2=q^{-1}a_1b_1,\quad d_{-2}=qa_{-1}b_{-1}.
\label{88}
\end{equation}
The choices \eqref{69} for $h_k,x_k,g_k$ corresponding to the
Askey--Wilson polynomials have the form \eqref{87} and satisfy the constraints
\eqref{88}. More generally,
Verde-Star \cite{1} claims that all families in the $q$-Askey scheme
\cite[Chapter 14]{2}, except for the continuous $q$-Hermite polynomials,
can be obtained in this way. This we have indeed shown in \cite{25}.

The eleven parameters $a_{-1},a_0,a_1$, $b_{-1},b_0,b_1$, $d_{-2},d_{-1},d_0,d_1,d_2$
with the three constraints $\sum_{i=-2}^2 d_i=0$ and \eqref{88}
can be further reduced
to four essential parameters by the four symmetries mentioned in \S\ref{91}.
Indeed, by the two translation symmetries the choices of $a_0$ and $b_0$ do not
really matter. The two dilation symmetries mean that the $a_i$ and $d_i$ may
be multiplied by the same nonzero factor, and that also the $b_i$ and $d_i$ may
be multiplied by the same nonzero factor.

There is a $q\leftrightarrow q^{-1}$ duality, in which
$a_j,b_j,d_j$ go to $a_{-j},b_{-j},d_{-j}$, respectively.
Finally the $x\leftrightarrow h$ duality \eqref{83} corresponds to an exchange
$a_i\leftrightarrow b_i$.

\subsection{The $q=1$ case}
Put $\xi=1+q+q^{-1}$ with $q=1$, so $\xi=3$.
For this case see \cite[\S8]{1}.
The solutions of \eqref{75}, \eqref{76}, \eqref{77} take the form
\begin{equation}
\label{89}
\begin{split}
&h_k=a_0+a_1 k+a_2 k^2,\qquad
x_k=b_0+b_1 k+b_2 k^2,\\
&g_k=d_1 k +d_2 k^2+d_3 k^3+d_4 k^4,\\
&a_1\notin a_2\ZZ_{<0},\mbox{ in particular,\quad $a_1$ or $a_2\ne 0$},
\quad \mbox{$d_i\ne0$ for some $i$}.
\end{split}
\end{equation}
Furthermore, the constraints \eqref{81}, \eqref{82} now simplify to
\begin{equation}
d_3=a_1b_2+a_2b_1-2a_2b_2,\quad d_4=a_2b_2.
\label{90}
\end{equation}
The choices \eqref{72} for $h_k,x_k,g_k$ corresponding to the
Wilson polynomials have the form \eqref{89} and satisfy the constraints
\eqref{90}. More generally,
Verde-Star \cite{1} claims that all families in the Askey scheme
\cite[Chapter 9]{2}, except for the Hermite polynomials,
can be obtained in this way. We will make this concrete in the next section.

The eight parameters $a_0,a_1,a_2$, $b_0,b_1,b_2$, $d_1,d_2$ can be reduced
to four essential parameters by the four symmetries mentioned in \S\ref{91}.
Indeed, by the two translation symmetries the choices of $a_0$ and $b_0$ do not
really matter. The two dilation symmetries mean that the $a_i$ and $d_i$ may
be multiplied by the same nonzero factor, and that also the $b_i$ and $d_i$ may
be multiplied by the same nonzero factor.
Note also that the $x\leftrightarrow h$ duality \eqref{83} corresponds to an exchange
$a_i\leftrightarrow b_i$.
\section{The Verde-Star scheme}
\subsection{Scheme with degree triples}
Recall \eqref{62}--\eqref{64}:
\begin{equation}
u_n(x)=\sum_{k=0}^n c_{n,k}\,v_k(x),\quad
v_k(x)=\prod_{j=0}^{k-1}(x-x_j),\quad
c_{n,k}=\prod_{j=k}^{n-1}\frac{g_{j+1}}{h_n-h_j}\,,\label{16}
\end{equation}
where we assume that $h_k,x_k,g_k$ are of the form \eqref{89} with constraints
\eqref{90}.

Associate with each given triple
$x_k,g_k,h_k$ of the form \eqref{89} the \emph{degree triple}
\begin{equation}
\deg(x),\deg(g),\deg(h),
\label{48}
\end{equation}
where $\deg(x)$ means the degree of $x_k$ as a polynomial in $k$, and similarly
for $\deg(g)$ and $\deg(h)$.
These three numbers are also the highest $i$ for which $b_i,d_i,a_i$, respectively,
are nonzero. We will classify the possible degree triples \eqref{48} and
associate them with families in the Askey scheme.
The degree triples have to satisfy the following rules, evident from
\eqref{89}, \eqref{90}.
\begin{enumerate}
\item
$\deg(x)\in\{0,1,2\}$, $\deg(g)\in\{1,2,3,4\}$, $\deg(h)\in\{1,2\}$.
\item
$\deg(g)=\deg(x)+\deg(h)$ if $\deg(x)+\deg(h)\ge3$,\\
$\deg(g)\le 2$ if $\deg(x)+\deg(h)\le2$.
\item
If for a given degree triple one of the degrees is decreased by 1 or, when
$\deg(x)+\deg(h)\ge3$, $\deg(g)$ and one of $\deg(x)$ and $\deg(h)$ are decreased
by 1, then we draw an arrow from that degree triple to the new one. It corresponds
to the situation that a $b_i$, $d_i$ or $a_i$, for highest~$i$ where it is
nonvanishing, is put to 0.
\item
If for a given degree triple the order is reversed then this means for corresponding
polynomials $u_n(x)$ that we are passing to the dual $\wt u_n(x)$, see
\eqref{83}.
\end{enumerate}

Figure \ref{21} gives the array of degree triples following the above rules.
To each degree triple in the scheme we have added one or more acronyms of
corresponding families of orthogonal polynomials in the Askey scheme:
\sLP
W (Wilson), R (Racah), cdH (continuous dual Hahn), dH (dual Hahn),
cH (continuous Hahn), H (Hahn), M--P (Meixner--Pollaczek), K (Krawtchouk),
J (Jacobi), Ch (Charlier), L (Laguerre), B (Bessel), and the degenerate case
bin (binomial, i.e., $(1-x)^n$, not in the Askey scheme).
\sLP
This perfectly matches with the usual Askey scheme \cite[Chapter 9]{2} except
for the Hermite polynomials, which do not fit into this framework.
The precise data $x_k,g_k,h_k$ for each family are given in the Appendix.
\setlength{\unitlength}{3mm}
\begin{figure}[h]
\centering
\begin{picture}(50,40)
\put(20,36.5)
{\framebox(4.5,3.5){\begin{tabular}{c}W, R\\2,4,2\end{tabular}}}
\put(20.1,36.4) {\vector(-2,-1){3.2}}
\put(10,31.3)
{\framebox(7,3.5){\begin{tabular}{c}cdH, dH\\2,3,1\end{tabular}}}
\put(24.6,36.4) {\vector(2,-1){3.2}}
\put(27.6,31.3)
{\framebox(5,3.5){\begin{tabular}{c}cH, H\\1,3,2\end{tabular}}}
\put(16.9,31.2) {\vector(1,-1){4.2}}
\put(27.6,31.2) {\vector(-1,-1){4.2}}
\put(18.8,23.6)
{\framebox(6.9,3.5){\begin{tabular}{c}M--P, M, K\\1,2,1\end{tabular}}}
\put(29.9,31.2) {\vector(1,-2){2}}
\put(29.4,23.6)
{\framebox(4.5,3.5){\begin{tabular}{c}J\\0,2,2\end{tabular}}}
\put(22,23.5) {\vector(0,-1){4}}
\put(25.6,23.5) {\vector(1,-1){4.2}}
\put(31.4,23.5) {\vector(0,-1){4}}
\put(20,15.9)
{\framebox(4.5,3.5){\begin{tabular}{c}Ch\\1,1,1\end{tabular}}}
\put(29.4,15.9)
{\framebox(4.5,3.5){\begin{tabular}{c}L\\0,2,1\end{tabular}}}
\put(34.0,23.5) {\vector(1,-1){4.2}}
\put(38.3,15.9)
{\framebox(4.5,3.5){\begin{tabular}{c}B\\0,1,2\end{tabular}}}
\put(24.6,15.8) {\vector(1,-1){4.2}}
\put(31.4,15.8){\vector(0,-1){4}}
\put(38.2,15.8) {\vector(-1,-1){4.2}}
\put(28.7,8.2)
{\framebox(6,3.5){\begin{tabular}{c}bin\\0,1,1\end{tabular}}}
\end{picture}
\vspace*{-2cm}
\caption{The Verde-Star scheme using degree triples}
\label{21}
\end{figure}
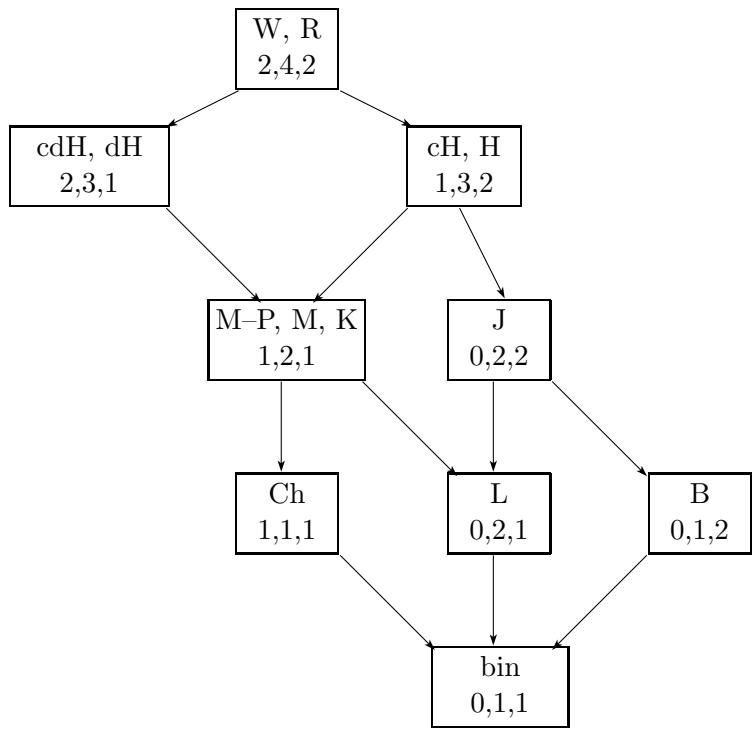

Note that the boxes for Wilson and Racah, for Meixner--Pollaczek, Meixner and
Krawtchouk, and for Charlier are self-dual, and that the boxes
for (continuous) dual Hahn and for (continuous) Hahn are dual to each other.
\subsection{Scheme with uniform parameters}
We have seen that a system of polynomials in the $q=1$ scheme is determined
by $h_k,x_k,g_k$ of the form \eqref{89} with constraints \eqref{90}.
Furthemore, putting $a_0=b_0=0$ does not change the system in an essential way.
So we are left with six parameters $a_1,a_2,b_1,b_2,d_1,d_2$ satisfying two
dilation symmetries: $a_1,a_2,d_1,d_2$ may be multiplied with the same factor and
$b_1,b_2,d_1,d_2$ may be multiplied with the same factor. Since it is not allowed 
that, for all $k$, $h_k$ is constant or $g_k=0$, we may assume that
$a_1\ne0$ and $d_1\ne0$. So we
can put $a_1=d_1=1$ by these dilation symmetries. Hence each system in the scheme
can be described by the four parameters $a_2,b_1,b_2,d_2$. 
From the rules we earlier gave for the degree triples it follows that
\begin{equation*}
d_2=0\;\Rightarrow\;a_2=b_2=0\;\;\mbox{or}\;\;a_1=a_2=0.
\end{equation*}
Also we will not let $a_1=0$ if $a_2\ne0$.

It turns out that each node in the scheme in Figure \ref{21} corresponds to
some of the four parameters being zero and that all allowed parameter restrictions
occur. We draw this scheme once more in Figure \ref{94}, but now in each box
with indication of the vanishing parameters.

\setlength{\unitlength}{3mm}
\begin{figure}[h]
\centering
\begin{picture}(50,40)
\put(20,36.5)
{\framebox{\begin{tabular}{c}W, R\end{tabular}}}
\put(20.1,35.6) {\vector(-2,-1){4}}
\put(10,31.3)
{\framebox{\begin{tabular}{c}cdH, dH\\$a_2=0$\end{tabular}}}
\put(25.3,35.6) {\vector(2,-1){4}}
\put(27.6,31.3)
{\framebox{\begin{tabular}{c}cH, H\\$b_2=0$\end{tabular}}}
\put(16.9,30) {\vector(1,-1){4.2}}
\put(27.6,30) {\vector(-1,-1){4.2}}
\put(18.8,23.6)
{\framebox{\begin{tabular}{c}M--P, M, K\\$a_2=b_2=0$\end{tabular}}}
\put(29.9,29.6) {\vector(1,-1){3.8}}
\put(29.4,23.6)
{\framebox{\begin{tabular}{c}J\\$b_1=b_2=0$\end{tabular}}}
\put(22,22) {\vector(-1,-1){4}}
\put(25.6,22) {\vector(1,-1){4}}
\put(31.4,22.1) {\vector(0,-1){4}}
\put(12,15.9)
{\framebox{\begin{tabular}{c}Ch\\$a_2=b_2=d_2=0$\end{tabular}}}
\put(25.4,15.9)
{\framebox{\begin{tabular}{c}L\\$a_2=b_1=b_2=0$\end{tabular}}}
\put(35.0,22) {\vector(1,-1){4}}
\put(38.3,15.9)
{\framebox{\begin{tabular}{c}B\\$b_1=b_2=d_2=0$\end{tabular}}}
\put(23.6,14.3) {\vector(1,-1){4}}
\put(31.4,14.3){\vector(0,-1){4}}
\put(38.2,14.3) {\vector(-1,-1){4}}
\put(24.7,8.2)
{\framebox{\begin{tabular}{c}bin\\$a_1=a_2=b_2=d_2=0$\end{tabular}}}
\end{picture}
\vspace*{-2cm}
\caption{The Verde-Star scheme with indication of vanishing parameters}
\label{94}
\end{figure}
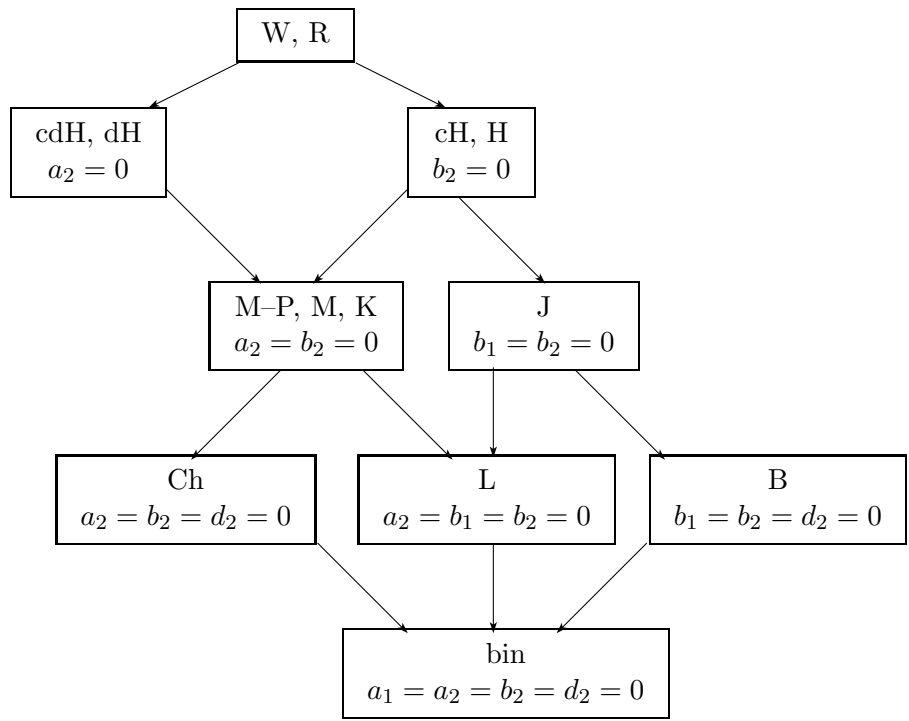

In order to get hypergeometric expressions for the families in the boxes of
Figure \ref{94} we have from \eqref{84}, \eqref{89}, \eqref{90} that
\begin{equation}
U_n(x;a_1,b_1,b_2,d_1)=
\sum_{k=0}^n\frac{(-n)_k}{k!}\,
\prod_{j=0}^{k-1}
\frac{\big(a_2(n+j)+1\big)(b_2j^2+b_1j-x)}
{a_2b_2(j+1)^3+(b_2+a_2b_1-2a_2b_2)(j+1)^2+d_2(j+1)+1}\,.
\label{96}
\end{equation}
Then put some parameters equal to zero according to what is indicated in a box
in Figure \ref{94}. However, rewriting \eqref{96} in the generic case
as a Wilson polynomial
with parameters $a,b,c,d$ expressed in terms of $a_1,b_1,b_2,d_1$ would be tedious.

We could also say that we have realized the Askey scheme (except for the Hermite
polynomials) as a very simple four-manifold (like $\RR^4$ or $\CC^4$) with
lower-dimensional submanifolds obtained by putting coordinates to zero.
Compare with an alternative and very complicated
description by the author in \cite{8}.
\section{Limit from the $q$-Verde-Star scheme to the Verde-Star scheme}
In the general $q$-case with $q\ne0$ and $1\notin q^\ZZ$
the $h_k,x_k,g_k$ are of general form \eqref{87}
with constraints \eqref{88}.
In \cite[Figure 1]{25} the (much larger) $q$-analogue of Figure \ref{21} is given
with arrays of black and white balls on the nodes. These balls correspond to
the parameters $b_i,d_i,a_i$ in \eqref{87}, white meaning a vanishing
parameter. Moreover most nodes in that diagram have a $q\leftrightarrow q^{-1}$ dual
giving a different array which is not drawn in \cite{25}.
Instead of working with black and white balls, we will use here a $q$-analogue of the
degree triples in Figure 1. Write \eqref{87} as
\begin{equation}
x_k=\sum_{j=\deg_-(x)}^{\deg_+(x)} b_j q^{kj},\quad
g_k=\sum_{j=\deg_-(g)}^{\deg_+(g)} d_j q^{kj},\quad
h_k=\sum_{j=\deg_-(h)}^{\deg_+(h)} a_j q^{kj},
\label{51}
\end{equation}
with $\deg_+(x), \deg_+(h)\in\{0,1\}$, $\deg_-(x), \deg_-(h)\in\{-1,0\}$,
$\deg_-(h)<\deg_+(h)$, and
$-2\le\deg_-(g)<\deg_+(g)\le2$, and such that $b_j$, $d_j$ and $a_j\ne0$
if $j=\deg_{\pm}(x)\ne0$, $j=\deg_{\pm}(g)$, and $j=\deg_{\pm}(h)\ne0$, respectively.
With $x_k,g_k,h_k$ as in \eqref{51} we associate a \emph{bidegree triple}
\begin{equation}
\deg_-(x), \deg_+(x); \deg_-(g), \deg_+(g); \deg_-(h),\deg_+(h)
\label{53}
\end{equation}
The black balls in the arrays in \cite[Figure 1]{25} run
from $\deg_-(x)$ to $\deg_+(x)$ in the first row,
from $\deg_-(g)$ to $\deg_+(g)$ in the second row,
and from $\deg_-(h)$ to $\deg_+(h)$ in the third row of the array.
In the $q\leftrightarrow q^{-1}$ duality
$b_j,d_j,a_j$ go to $b_{-j},d_{-j},a_{-j}$, respectively.
Accordingly there are interchanges $\deg_-(x)\leftrightarrow -\deg_+(x)$,
$\deg_-(g)\leftrightarrow -\deg_+(g)$, $\deg_-(h)\leftrightarrow -\deg_+(h)$.

Now we want to let $q\to 1$. Then, in many cases but not always,
the bidegree triple \eqref{53} goes to
the degree triple \eqref{48} with
$\deg(x)=\deg_+(x)-\deg_-(x)$, $\deg(g)=\deg_+(g)-\deg_-(g)$, and
$\deg(h)=\deg_+(h)-\deg_-(h)$.

\begin{Example}[Askey--Wilson $\to$ Wilson]\rm
Consider the Askey--Wilson polynomial $U_n(x;a,b,c,d;q)$ in
\eqref{92} with data \eqref{69}. Then the rescaled polynomial
\begin{align*}
&U_n\big(2-(1-q)^2 x;q^a,q^b,q^c,q^d;q\big)\\
&=\sum_{k=0}^n\,\frac{(1-q)^{2k}(q^{-n},q^{n+a+b+c+d-1};q)_k}
{(q^{a+b},q^{a+c},q^{a+d},q;q)_k}\,
q^{k(a+1)+\half k(k-1)}
\prod_{j=0}^{k-1}
\left(x+\frac{(q^{\half(a+j)}-q^{-\half(a+j)})^2}{(1-q)^2}\right)
\end{align*}
has corresponding data
\begin{align*}
&h_k=-\frac{q^{-k}(1-q^k)(1-q^{k+a+b+c+d-1})}{(1-q)^2}\,,\quad
x_k=-\frac{\big(q^{\half(a+k)}-q^{-\half(a+k)}\big)^2}{(1-q)^2}\,,\\
&g_k=\frac{q^{-a-2k+1}(1-q^{k+a+b-1})(1-q^{k+a+c-1})(1-q^{k+a+d-1})(1-q^k)}{(1-q)^4}\,,
\end{align*}
where we have multiplied $h_k$ and $g_k$ with the same suitable number.
On letting $q\to 1$ we arrive at the data \eqref{72} for Wilson polynomials.
Hence
\begin{equation}
\lim_{q\to1} U_n\big(2-(1-q)^2 x;q^a,q^b,q^c,q^d;q\big)=U_n(x;a,b,c,d),
\end{equation}
where $U_n$ on the left-hand side is given by \eqref{92} and on the
right-hand side by \eqref{93}.

In the limit formula \cite[(14.1.21)]{2} from Askey--Wilson to Wilson there is a
non-algebraic rescaling of the $z$-variable as $q^{\iup x}$. 
Similar non-algebrtaic rescalings occur on many other places in
\cite[Chapter 14]{2}, but in all such cases an algebraic rescaling will be possible
(see some examples in Appendix B).

Furthermore note that
there is no simple rescaling of $a_{-1},a_0,a_1$, $b_{-1},b_0,b_1$,
$d_{-2},d_{-1},d_0,d_1,d_2$ in \eqref{87} to enable the limit for $q\to1$.
This also means that there is no simultaneous parametrization of the $q$-case
\eqref{87} and the $q=1$ case \eqref{90} such that limits for $q\to1$ can be taken
in this parametrization.
\end{Example}

\subsection{A list of $q\to1$ limits}
In the following list we give, for each node in the Verde-Star scheme Figure \ref{21},
one or more nodes in the $q$-Verde-Star scheme \cite[Figure 1]{25} for which the data
$x_k,h_k,g_k$ tend to the data for the node in Figure \ref{21} as $q\to1$
after suitable
rescaling and such that the number of (essential) parameters is preserved in the
limit (i.e., the nodes have to be in the same row in the two figures).
The limit for the data will imply a limit for the corresponding polynomials.
In the $q=1$ case we also give the degree triples and in the $q$-case we also give
the bidegree triples (taken from the
arrays of black and white balls in \cite[Figure 1]{25}). One could have added the
bidegree triples obtained by $q\leftrightarrow q^{-1}$ exchange.
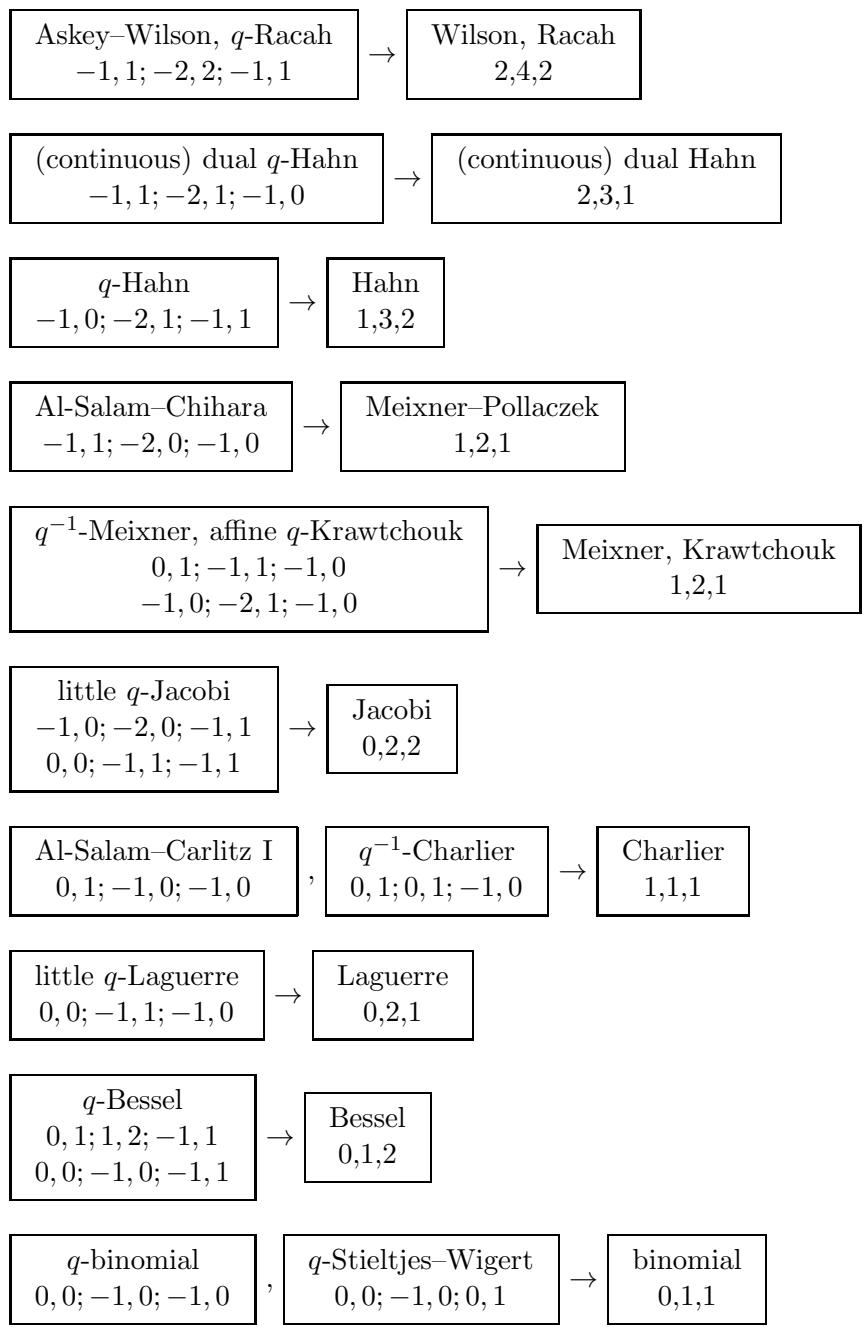
\bLP
\framebox{\begin{tabular}{c}Askey--Wilson, $q$-Racah\\$-1,1;-2,2;-1,1$\end{tabular}}
$\to$
\framebox{\begin{tabular}{c}Wilson, Racah\\2,4,2\end{tabular}}
\bLP
\framebox{\begin{tabular}{c}(continuous) dual $q$-Hahn\\
$-1,1;-2,1;-1,0$\end{tabular}}
$\to$
\framebox{\begin{tabular}{c}(continuous) dual Hahn\\2,3,1\end{tabular}}
\bLP
\framebox{\begin{tabular}{c}$q$-Hahn\\
$-1,0;-2,1;-1,1$\end{tabular}}
$\to$
\framebox{\begin{tabular}{c}Hahn\\1,3,2\end{tabular}}
\bLP
\framebox{\begin{tabular}{c}Al-Salam--Chihara\\
$-1,1;-2,0;-1,0$\end{tabular}}
$\to$
\framebox{\begin{tabular}{c}Meixner--Pollaczek\\1,2,1\end{tabular}}
\bLP
\framebox{\begin{tabular}{c}$q^{-1}$-Meixner, affine $q$-Krawtchouk\\
$0,1;-1,1;-1,0$\\$-1,0;-2,1;-1,0$\end{tabular}}
$\to$
\framebox{\begin{tabular}{c}Meixner, Krawtchouk\\1,2,1\end{tabular}}
\bLP
\framebox{\begin{tabular}{c}little $q$-Jacobi\\
$-1,0;-2,0;-1,1$\\$0,0;-1,1;-1,1$\end{tabular}}
$\to$
\framebox{\begin{tabular}{c}Jacobi\\0,2,2\end{tabular}}
\bLP
\framebox{\begin{tabular}{c}Al-Salam--Carlitz I\\$0,1;-1,0;-1,0$\end{tabular}} ,
\framebox{\begin{tabular}{c}$q^{-1}$-Charlier\\$0,1;0,1;-1,0$\end{tabular}}
$\to$
\framebox{\begin{tabular}{c}Charlier\\1,1,1\end{tabular}}
\bLP
\framebox{\begin{tabular}{c}little $q$-Laguerre\\$0,0;-1,1;-1,0$\end{tabular}}
$\to$
\framebox{\begin{tabular}{c}Laguerre\\0,2,1\end{tabular}}
\bLP
\framebox{\begin{tabular}{c}$q$-Bessel\\$0,1;1,2;-1,1$\\$0,0;-1,0;-1,1$\end{tabular}}
$\to$
\framebox{\begin{tabular}{c}Bessel\\0,1,2\end{tabular}}
\bLP
\framebox{\begin{tabular}{c}$q$-binomial\\$0,0;-1,0;-1,0$\end{tabular}} ,
\framebox{\begin{tabular}{c}$q$-Stieltjes--Wigert\\$0,0;-1,0;0,1$\end{tabular}}
$\to$
\framebox{\begin{tabular}{c}binomial\\0,1,1\end{tabular}}
\bPP
Some worked-out examples are given in Appendix B.
\appendix
\section{Explicit data for the families in Figure \ref{21}}
\label{61}
For each family we start with the acronym, followed by the corresponding degree
triple. After the full name of the family the corresponding section number in
\cite[Chapter 9]{2} is given.
\bLP
\textbf{W(2,4,2) Wilson} (1): $u_n(x)=k_n^{-1} W_n(x;a,b,c,d)$,\\
$x_k=-(k+a)^2$,\quad
$g_k=k(k+a+b-1)(k+a+c-1)(k+a+d-1)$,\\
$h_k=-k(k+a+b+c+d-1)$,\quad
$k_n=(-1)^n (n+a+b+c+d-1)_n$.
\mLP
\textbf{R(2,4,2) Racah} (2): $u_n(x)=k_n^{-1} R_n(x;\al,\be,\ga,\de)$\quad
($\al+1$ or $\be+\de+1$ or $\ga+1=-N$),\\
$x_k=k(k+\ga+\de+1)$,\quad
$g_k=k(k+\al)(k+\be+\de)(k+\ga)$,\quad
$h_k=k(k+\al+\be+1)$,\sLP
$k_n=\dstyle\frac{(n+\al+\be+1)_n}{(\al+1)_n (\be+\de+1)_n (\ga+1)_n}$\,.
\mLP
\textbf{cdH(2,3,1) continuous dual Hahn} (3): $u_n(x)=k_n^{-1} S_n(x;a,b,c)$,\\
$x_k=-(k+a)^2$,\quad
$g_k=k(k+a+b-1)(k+a+c-1)$,\quad
$h_k=-k$,\quad
$k_n=(-1)^n$.
\mLP
\textbf{dH(2,3,1) dual Hahn} (6): $u_n(x)=k_n^{-1} R_n(x;\ga,\de,N)$,\\
$x_k=k(k+\ga+\de+1)$,\quad
$g_k=k(k+\ga)(N-k+1)$,\quad
$h_k=-k$,\quad
$k_n=\dstyle\frac1{(\ga+1)_n(-N)_n}$\,.
\mLP
\textbf{cH(1,3,2) continuous Hahn} (4): $u_n(x)=k_n^{-1} p_n(x;a,b,c,d)$\quad
($c=\overline a$ and $d=\overline b$),\\
$x_k=\iup(k+a)$,\quad
$g_k=\iup k(k+a+c-1)(k+a+d-1)$,\quad
$h_k=k(k+a+b+c+d-1)$,\sLP
$k_n=\dstyle\frac{(n+a+b+c+d-1)_n}{n!}$\,.
\mLP
\textbf{H(1,3,2) Hahn} (5): $u_n(x)=k_n^{-1} Q_n(x;\al,\be,N)$,\\
$x_k=k$,\quad
$g_k=k(k+\al)(k-N-1)$,\quad
$h_k=k(k+\al+\be+1)$,\quad
$k_n=\dstyle\frac{(n+\al+\be+1)_n}{(\al+1)_n (-N)_n}$\,.
\mLP
\textbf{M--P(1,2,1) Meixner--Pollaczek} (7): $u_n(x)=k_n^{-1} P_n^{(\la)}(x;\phi)$,\\
$x_k=\iup(k+\la)$,\quad
$g_k=k(k+2\la-1)$,\quad
$h_k=-\iup(1-\eup^{-2\iup\phi})(k+\la)$,\quad
$k_n=\dstyle\frac{2^n(\sin\phi)^n}{n!}$\,.
\mLP
\textbf{M(1,2,1) Meixner} (10): $u_n(x)=k_n^{-1} M_n(x;\be,c)$,\\
$x_k=k$,\quad
$g_k=k(k+\be-1)$,\quad
$h_k=(1-c^{-1})k$,\quad
$k_n=\dstyle\frac{(1-c^{-1})^n}{(\be)_n}$\,.
\mLP
\textbf{K(1,2,1) Krawtchouk} (11): $u_n(x)=k_n^{-1} K_n(x;p,N)$,\sLP
$x_k=k$,\quad
$g_k=k(k-N-1)$,\quad
$h_k=p^{-1}k$,\quad
$k_n=\dstyle\frac1{p^n (-N)_n}$\,.
\mLP
\textbf{J(0,2,2) Jacobi} (8): $u_n(x)=k_n^{-1} P_n^{(\al,\be)}(x)$,\\
$x_k=1$,\quad
$g_k=k(k+\al)$,\quad
$h_k=\thalf k(k+\al+\be+1)$,\quad
$k_n=\dstyle\frac{(n+\al+\be+1)_n}{2^n n!}$\,.
\mLP
\textbf{Ch(1,1,1) Charlier} (14): $u_n(x)=k_n^{-1} C_n(x;a)$,\\
$x_k=k$,\quad
$g_k=k$,\quad
$h_k=-a^{-1}k$,\quad
$k_n=(-a)^{-n}$\,.
\mLP
\textbf{L(0,2,1) Laguerre} (12): $u_n(x)=k_n^{-1} L_n^{(\al)}(x)$,\sLP
$x_k=0$,\quad
$g_k=k(k+\al)$,\quad
$h_k=-k$,\quad
$k_n=\dstyle\frac{(-1)^n}{n!}$\,.
\mLP
\textbf{B(0,1,2) Bessel} (13): $u_n(x)=k_n^{-1} y_n(x;a)$,\\
$x_k=0$,\quad
$g_k=k$,\quad
$h_k=\thalf k(k+a+1)$,\quad
$k_n=\dstyle\frac{(n+a+1)_n}{2^n}$\,.
\mLP
\textbf{bin(0,1,1) binomial} (degenerate): $u_n(x)=k_n^{-1} (1-x)^n$,\\
$x_k=0$,\quad
$g_k=k$,\quad
$h_k=-k$,\quad
$k_n=(-1)^n$\,.
\section{Some worked-out $q\to1$ limits}
In each case we start in the $q$-case with $U_n(x)$, given by \eqref{84} and
written as $q$-hypergeometric function, and with the corresponding data taken from
\cite[Appendix A]{25}. Then we rescale $U_n(x)$, write the data corresponding
to the rescaling (and perform on $h_k$ a possible translation and on $h_k$ and $g_k$
a multiplication with the same factor), let $q\to1$, and obtain the data for the
limit polynomials which turn out to be equal to the data given in Appendix A.

\paragraph{Al-Salam--Chihara $\to$ Meixner--Pollaczek} \cite[(14.8.19)]{2}
\begin{align*}
&U_n(x;a,b;q)=\qhyp32{q^{-n},az,az^{-1}}{ab,0}{q,q}
=Q_n(\thalf x;a,b\,|\,q)\quad(x=z+z^{-1}),\\
&h_k=q^{-k}-1,\quad
x_k=aq^k+a^{-1}q^{-k},\quad
g_k=q^{-2k+1}a^{-1}(1-abq^{k-1})(1-q^k).\sLP
&U_n\big(2\cos\phi+2x(1-q)\sin\phi;
q^\la\eup^{-\iup\phi},q^\la\eup^{\iup\phi};q\big)=\sum_{k=0}^n
\frac{(q^{-n};q)_k(1-q)^k}{(q^{2\la},q;q)_k}\,\big(\iup(1-\eup^{-2\iup\phi})\big)^k\,
\prod_{j=0}^{k-1}(x-x_j),\\
&h_k=-\frac{1-q^{-k}}{1-q},\quad
x_k=-\iup\,\frac{1-q^{-\la-k}}{1-q}\,
\frac{1-q^{\la+k}\eup^{-2\iup\phi}}{1-\eup^{-2\iup\phi}}\,,\quad
g_k=\frac{q^{-\la-2k-1}}{2\eup^{-\iup\phi}\sin\phi}\,
\frac{(1-q^{2\la+k-1})(1-q^k)}{(1-q)^2}\,.\sLP
&U_n(x;\la,\phi)=\hyp21{-n,\la+\iup x}{2\la}{1-\eup^{-2\iup\phi}}=
\frac{\eup^{-\iup n\phi} n!}{(2\la)_n}\,P_n^{(\la)}(x;\phi),\\
&h_k=k,\quad
x_k=\iup(\la+k),\quad
g_k=\frac{k(2\la+k-1)}{2\eup^{-\iup\phi}\sin\phi}\,.
\end{align*}

\paragraph{Al-Salam--Carlitz $\to$ Charlier} \cite[(14.24.13)]{2}
\begin{align*}
&U_n(x;a;q)=\qhyp21{q^{-n},x^{-1}}0{q,qa^{-1}x}
=(-a)^{-n} q^{-\half n(n-1)} U_n^{(a)}(x;q),\\
&h_k=-a^{-1}(1-q^{-k}),\quad
x_k=q^k,\quad
g_k=1-q^{-k}.\sLP
&U_n\big(1-(1-q)x;-(1-q)a;q\big)=\sum_{k=0}^n
\frac{(q^{-n};q)_k}{(q;q)_k}\,(-qa^{-1})^k\,
\prod_{j=0}^{k-1}\left(-x+\frac{1-q^j}{1-q}\right),\\
&h_k=a^{-1}\,\frac{1-q^{-k}}{1-q},\quad
x_k=\frac{1-q^k}{1-q},\quad
g_k=-\frac{1-q^{-k}}{1-q}\,.\sLP
&U_n(x;a)=\hyp20{-n,-x}-{-a^{-1}}=C_n(x;a),\quad
h_k=-a^{-1}k,\quad
x_k=k,\quad
g_k=k.
\end{align*}

\paragraph{Stieltjes--Wigert $\to$ binomial}
\begin{align*}
&U_n(x;q)
=\qhyp11{q^{-n}}0{q,-q^{n+1}x}=S_n(x;q),\quad
x_k=0,\quad
h_k=-\frac{1-q^k}{1-q}\,,\quad
g_k=-\frac{1-q^{-k}}{1-q}\,.\\
&U_n(x)
=\hyp10{-n}-x=(1-x)^n,\quad
x_k=0,\quad
h_k=-k,\quad
g_k=k.
\end{align*}

\quad\\
\begin{footnotesize}
\begin{quote}
{ T. H. Koornwinder, Korteweg-de Vries Institute, University of
 Amsterdam,\\
 P.O.\ Box 94248, 1090 GE Amsterdam, The Netherlands;

\vspace{\smallskipamount}
email: }\url{thkmath@xs4all.nl}
\end{quote}
\end{footnotesize}

\end{document}